\documentclass[preprint,12pt]{elsarticle}

\usepackage{amssymb}
\usepackage{amsfonts}
\usepackage{amsmath}
\usepackage{float}
\usepackage{color}
\usepackage{url}
\usepackage{epstopdf}
\newtheorem{lemma}{Lemma}
\newtheorem{remark}{Remark}
\newtheorem{theorem}{Theorem}
\newtheorem{proposition}{Proposition}
\newtheorem{corollary}{Corollary}

\newcommand{\Z}{ \mathbb{Z} }

\usepackage{tcolorbox}
\usepackage{epstopdf}
\usepackage{subcaption}
\journal{Statistics and Probability Letters}

\begin{document}

\begin{frontmatter}



\title{A note on strong-consistency of   componentwise ARH(1) predictors}


\author{M. D. Ruiz-Medina \dag and J. \'Alvarez-Li\'ebana}

\address{\dag Department of Statistics and Operation
Research (mruiz@ugr.es,  javialvaliebana@ugr.es)\\  Faculty of Sciences,  University of Granada\\
Campus Fuente Nueva s/n\\18071 Granada, Spain}


\begin{abstract}
This paper presents a  new result on strong-consistency, in the trace norm,  of a  diagonal componentwise parameter estimator of the autocorrelation o-\linebreak perator of an autoregressive process of order one  (ARH(1) process), allowing strong-consistency of the  associated plug-in predictor.  These results are derived, when the  eigenvectors of the autocovariance operator are unknown, and  the autocorrelation operator does not admit a diagonal spectral representation with respect to the eigenvectors of the autocovariance operator.     
\end{abstract}

\begin{keyword}
Dimension reduction techniques  \sep empirical orthogonal bases \sep functional prediction \sep   
 strong-consistency \sep trace norm

\vspace*{0.5cm}



\noindent \textit{2010 Mathematics Subject Classification}
Primary 60G10; 60G15. Secondary 60F99; 60J05; 65F15
\end{keyword}

\end{frontmatter}


\section{Introduction.}
\label{sec:1}

 It is well-known  that functional prediction from the ARH(1) framework  (see Bosq \cite{Bosq00} and   Bosq and  Blanke \cite{Bosq07}) has become a very active  research area, driven by its interesting applications in the analysis of high-dimensional data  (see, for example, \'Alvarez-Li\'ebana, Bosq and  Ruiz-Medina  \cite{Alvarez1600}). Several authors have studied the asymptotic properties of componentwise estimators of the autocorrelation operator, and ARH(1) predictors, in the case of known and unknown eigenvectors of the autocovariance operator. We refer to \cite{Guillas01,Mas99,Mas04,Mas07}, where the efficiency, consistency and asymptotic normality of these estimators are addressed, in a parametric framework (see also \'Alvarez-Li\'ebana, Bosq and  Ruiz-Medina  \cite{Alvarez16}, on weak consistency in the norm of  Hilbert-Schmidt operators). 
   In this paper, we pay attention to  the diagonal componentwise consistent estimation of the autocorrelation operator of an ARH(1) process, from  projection into the  empirical eigenvectors of the autocovariance operator, when their theoretical counterparts do not diagonalize the autocorrelation operator.
   A suitable
sample-size-dependent truncation order, according to the spectral properties of the autocovariance operator of an ARH(1) process (namely, decay velocity and distance between eigenvalues), is usually selected to obtain the desirable asymptotic properties of componentwise   parameter  estimators, and their associated plug-in predictors, under different setting of conditions. The usual assumption on the autocorrelation operator is its boundedness or Hilbert-Schmidt property. Hence, the error term is usually measured in both operator norms, the norm of the space of bounded linear operators and the Hilbert-Schmidt norm.   In the same way, dimension reduction has been a crucial topic in the nonparametric and semiparametric frameworks. In these contexts, a kernel-based predictor is usual adopted, as early formulated in \cite{Poggi94}. Additional covariates were incorporated in the novel semiparametric kernel-based proposal by \cite{Aneiros08}, in which an extension, to the functional time series framework, of partial linear models, was firstly developed. In both  frameworks, regular autocorrelation operators usually arise. Thus, the  error term could be measured in the trace norm, in terms of the diagonal values, with the subsequent  dimension reduction. This is the proposal of the note carried out here, which  could help the usual computational limitations arising in the implementation of kernel-based estimators  (see also \cite{Antoniadis06}). A fractal approach was proposed in \cite{Ferraty02} to solve the curse of dimensionality, by imposing a concentration assumption about the distribution of the stochastic process.  For example, in the Gaussian case,  under the trace class assumption on the autocorrelation operator,  we can go beyond the concentration assumption on the autocorrelation kernel values, characterizing the distributional properties of the underlying stochastic process. A diagonal design of the kernel-based estimator can then be considered, ensuring the desirable asymptotic properties.

This paper intends to focus the attention of the reader on  the smoothing and regularization techniques usually applied in the functional data context, that allow to work on a  trace-class autocorrelation operator context, and hence, to measure the autocorrelation operator  estimation error, in the trace norm.     
Particularly, strong-consistency in the trace norm also implies the consistency  in the norm of the space of Hilbert-Schmidt operators, as well as in the norm of bounded linear operators.
 Results derived in this paper are mainly motivated by this fact, allowing an important  dimension reduction, in the componentwise estimation of the autocorrelation operator of an ARH(1) process, and its plug-in prediction.    It is important to highlight that the distributional properties of a wide class of stochastic processes in time and/or space can be characterized in terms of stochastic evolution equations, involving trace and positive autocorrelation operators   (see, for example,  Anh, Leonenko and Ruiz-Medina \cite{Anh16}, in the spatiotemporal case).
\section{Preliminaries.}
\label{sec:2}

 Let  $H$ be a real separable Hilbert space, and let  $X = \left\lbrace X_n,\ n \in \mathbb{Z} \right\rbrace $
 be  a zero-mean ARH(1) process on the basic probability space $(\Omega,\mathcal{A},P)$ satisfying the following  equation:
\begin{equation}
X_{n } = \rho \left(X_{n-1} \right) + \varepsilon_{n},\ n \in
\Z,\label{24bb}
\end{equation}
\noindent  where $\rho $ is the  autocorrelation operator of
 process $X$, which belongs to the space $\mathcal{L}(H)$ of bounded
linear operators,
such that $\Vert \rho^{k}
\Vert_{\mathcal{L}\left(H\right)}<1,$  for $k\geq k_{0},$ and for some
 $k_{0}.$ Here, as usual,  $\|\cdot\|_{\mathcal{L}(H)}$  denotes
the norm in the space $\mathcal{L}(H)$. The Hilbert-valued
innovation process $\varepsilon= \left \lbrace \varepsilon_{n}, \ n\in \mathbb{Z}\right\rbrace$ is assumed to be a strong white noise, and to be  uncorrelated with the random
initial condition. That is,
$\varepsilon $ is a Hilbert-valued zero-mean stationary process,
with independent and identically distributed components in time, and
with \linebreak $\sigma^{2}_{\varepsilon}={\rm E} \left[
\|\varepsilon_{n}\|_{H}^{2} \right] <\infty,$  for all $n\in
\mathbb{Z}.$ We restrict our attention here  to the case where
$\rho$ is such that  $\|\rho\|_{\mathcal{L}(H)}<1$.

Under the above-setting of conditions, involved  in the introduction of equation (\ref{24bb}), $X$  admits the following MAH($\infty$) representation (see
Bosq \cite{Bosq00}):
\begin{equation}
X_n = \displaystyle \sum_{k=0}^{\infty} \rho^{k} \left( \varepsilon_{n-k} \right), \quad n \in \mathbb{Z}, \label{2}
\end{equation}
\noindent  that provides  the unique stationary solution to equation (\ref{24bb}).
The autocovariance operator $C_{X}$ of the ARH(1) process $X$ is given by $$C_{X}={\rm E}[X_{n}\otimes X_{n}]=
{\rm E}[X_{0}\otimes X_{0}],\quad \forall n\in \mathbb{Z},$$
\noindent which is assumed to be a trace operator.

 Let us  now consider the empirical autocovariance operator $C_n$ of the ARH(1) process $X=\{ X_{n},\ n\in \mathbb{Z}\},$  computed as follows: 
\begin{equation}
C_n = \frac{1}{n} \displaystyle \sum_{i=0}^{n-1} X_i \otimes X_i,\quad n \geq 2,  \label{66aa}
\end{equation}
\noindent from a functional sample $X_{0}, X_{1},\dots, X_{n-1},$ of size $n$, of the ARH(1) process $X.$  The empirical eigenvalues $\{C_{n,j}, \ j \geq 1\},$ and  eigenvectors $\{\phi_{n,j},\ j\geq 1\}$ of $C_{n}$ are  given by (see \cite{Bosq00}, pp. 102--103):
 \begin{eqnarray}
C_{n}\phi_{n,j}&=& C_{n,j}\phi_{n,j},\quad j\geq 1,\quad C_{n,1}\geq \dots\geq C_{n,n}\geq 0=C_{n,n+1}=C_{n,n+2}=\dots \nonumber\\ && \label{ev1}
\end{eqnarray}
\noindent
In the remainder, we will denote
 \begin{eqnarray}X_{i,j,n} &=& \langle X_i, \phi_{n,j} \rangle_H,\quad i\in \mathbb{Z},~j \geq 1,~ n\geq 2,\nonumber\\
 \phi_{n,j}^{\prime }&=& \mbox{sgn}\left\langle \phi_{n,j} , \phi_{j} \right\rangle_{H}\phi_{j},\quad j \geq 1,~ n\geq 2,
 \label{teig}
 \end{eqnarray}
 \noindent where $\mbox{sgn} \langle \phi_{n,j}, \phi_j \rangle_H= \mathbf{1}_{\langle \phi_{n,j}, \phi_j \rangle_H\geq 0}-\mathbf{1}_{\langle \phi_{n,j}, \phi_j \rangle_H< 0}.$
Thus, \begin{equation}C_{n,j} =\frac{1}{n} \displaystyle \sum_{i=0}^{n-1} X_{i,j,n}^{2},\quad j \geq 1,\ n\geq 2.\label{eive1}\end{equation}

Let us now consider  the cross-covariance operator $D_{X}$ given by 
$$D_{X}={\rm E}[X_{i}\otimes X_{i+1}]={\rm E}[X_{0}\otimes X_{1}],\quad \forall i\in\mathbb{Z},$$
\noindent which is assumed to be a trace operator. The empirical cross-covariance operator $D_n$
is defined as follows: 
\begin{eqnarray}
D_n &=& \frac{1}{n-1}\displaystyle \sum_{i=0}^{n-2} X_i \otimes X_{i+1}, \quad n \geq 2.\label{37aaa}\end{eqnarray}
 Given the empirical eigenvectors of $C_{n},$ $\left\lbrace \phi_{n,j},~j \geq 1 \right\rbrace,$
  we compute 
\begin{eqnarray}
& & D_{n,j,l}^{*} = \langle D_n \left(\phi_{n,j} \right), \phi_{n,l}  \rangle_H= \frac{1}{n-1}\displaystyle
\sum_{i=0}^{n-2}X_{i,j,n}X_{i+1,l,n},~j,l \geq 1, \ n\geq 2,\nonumber
\end{eqnarray}
\noindent and, in particular,  we will use the notation
\begin{equation}
D_{n,j}=D_{n,j,j}^{*} =  \langle D_n \left(\phi_{n,j} \right), \phi_{n,j}  \rangle_H = \frac{1}{n-1} \displaystyle \sum_{i=0}^{n-2} X_{i,j,n}X_{i+1,j,n},\quad j \geq 1,\ n \geq 2,\label{38aaa}
\end{equation}
\noindent where the empirical projections $\{X_{i,j,n}, \ j\geq 1,~i=0,1,\ldots,n-1\}$ have been previously defined in equation (\ref{teig}).
\section{Auxiliary results}
\label{ar}
The  auxiliary results presented  in this section allow the derivation of Theorem \ref{proposition5}  in the next section. 
\medskip

\noindent  \textbf{Assumption A1.} The   random initial condition  $X_{0}$ of the ARH(1) process defined in (\ref{24bb}) satisfies
$${\rm E} \left\lbrace \left\| X_0 \right\|_{H}^{4} \right\rbrace < \infty.$$

\begin{theorem}
\label{theorem2} (see Theorem 4.1 on pp. 98--99, Corollary 4.1 on pp. 100--101 and  Theorem 4.8 on pp. 116--117, in \cite{Bosq00}).
Under \textbf{Assumption A1}, for any $\beta > \frac{1}{2},$ as $n\rightarrow \infty,$
\begin{equation}
\frac{n^{1/4}}{\left(\ln(n) \right)^{\beta}}  \left\| C_n - C_{X} \right\|_{\mathcal{S}(H)} \to^{a.s.}0, \quad \frac{n^{1/4}}{\left(\ln(n) \right)^{\beta}} \left\| D_n - D_{X} \right\|_{\mathcal{S}(H)} \to^{a.s.}0,
\end{equation}
\noindent and, if $\left\| X_0 \right\|_{H}$ is bounded,
\begin{eqnarray}
& &\left\| C_n - C_{X} \right\|_{\mathcal{S}(H)} = \mathcal{O} \left(\left(\frac{\ln(n) }{n} \right)^{1/2} \right)~a.s.,
\nonumber\\
& &\left\| D_n - D_{X} \right\|_{\mathcal{S}(H)} = \mathcal{O} \left(\left(\frac{\ln(n) }{n} \right)^{1/2} \right)~a.s.,
\end{eqnarray}
\noindent where $\left\| \cdot \right\|_{\mathcal{S}(H)}$ denotes the norm of the  Hilbert-Schmidt operators  on $H,$ and a.s. means, as usual,  almost surely. \end{theorem}

In the following, we denote by $\left\lbrace C_j, \ j \geq 1 \right\rbrace$ the sequence of eigenvalues of the autocovariance operator $C_{X}$, satisfying
\begin{equation}
C_{X}(\phi_{j})=C_{j}\phi_{j},\quad j\geq 1, \label{eq_new}
\end{equation}
\noindent with $\{\phi_{j},\ j\geq 1\}$ being the associated system of eigenvectors. 
\begin{lemma}
\label{lem1}
Under \textbf{Assumption A1}, for $n$ sufficiently large,
\begin{eqnarray}
\frac{n^{1/4}}{\left(\ln(n) \right)^{\beta}} \displaystyle \sup_{j \geq 1} \left| C_{n,j}- C_j \right|  & \leq & \frac{n^{1/4}}{\left(\ln(n) \right)^{\beta}} \left\| C_n - C \right\|_{\mathcal{S}(H)} \to^{a.s.} 0,\label{c11b}\end{eqnarray}
\noindent where $\{C_{j}, \ j\geq 1\}$ have been introduced in equation (\ref{eq_new}). 
\end{lemma}

The proof of Lemma \ref{lem1} is straightforward since, from Theorem \ref{theorem2} and $n$ sufficiently large,  $C_{n}$  is a  Hilbert-Schmidt operator, and in particular, it is a compact operator. Thus, applying Lemma 4.2 on p. 103 in \cite{Bosq00}, and Theorem \ref{theorem2},  for $n\geq n_{0},$ with $n_{0}$ sufficiently large, we obtain
\begin{eqnarray}
 \frac{n^{1/4}}{\left(\ln(n) \right)^{\beta}}\sup_{k\geq 1}|C_{n,k}-C_{k}|
 &\leq & \frac{n^{1/4}}{\left(\ln(n) \right)^{\beta}}\|C_{n}-C_{X}\|_{\mathcal{S}(H)}\to^{a.s.} 0.
\label{eqdpc1}
\end{eqnarray}

Proposition \ref{prr1} below obtains  the strong-consistency  of $\left\{D_{n,j}, \ j \geq 1\right\},$ introduced in equation (\ref{38aaa}). In the derivation of these strong-consistency results, the distance between the  eigenvalues of the autocovariance operator, and their rate of convergence to zero, will play a key role. Namely, if the following quantities are considered:
\begin{equation}\Lambda_{k}=\sup_{1\leq j\leq k}(C_{j}-C_{j+1})^{-1},\quad k\geq 1,\label{uee}
\end{equation}
\noindent the rate of divergence of $\Lambda_{k_n}$ is crucial, as can be noted in the statements of Lemma \ref{lemmanew} and Proposition \ref{prr1} (see also Remark \ref{remratec} below). Henceforth, $k_n$ denotes a truncation parameter, which verifies
\begin{equation}
\displaystyle \lim_{n \to \infty} k_n = \infty, \quad \frac{k_n}{n} < 1, \quad k_n \geq 1. \label{trunc}
\end{equation}

\begin{lemma}
\label{lemmanew}
 Let us consider $\Lambda_{k_{n}}$ introduced in equation (\ref{uee}), for a given truncation parameter $k_n$, as reflected in (\ref{trunc}).
Assuming that $\| X_{0}\|_{H}$ is bounded and
 $\Lambda_{k_{n}}=o\left(n^{1/4}(\ln(n))^{\beta -1/2} \right),$ as $n\rightarrow \infty$, under \textbf{Assumption A1}, the following limit then holds
\begin{equation}\frac{n^{1/4}}{\left(\ln(n) \right)^{\beta}}\sup_{1\leq j\leq k_{n}}\|\phi_{n,j}^{\prime }-\phi_{n,j}\|_{H}\to^{a.s.} 0,\quad n\rightarrow \infty,\label{leqkk}\end{equation}
 \noindent for any $\beta > 1/2$, where $\{\phi_{n,j}^{\prime },\ j\geq 1\}$ are  introduced in equation (\ref{teig}).
\end{lemma}

\begin{proof}
From equation (4.44) in Lemma 4.3 on page 104 of  \cite{Bosq00}, for any $n \geq 2$ and $1 \leq j \leq k_n$,
  \begin{eqnarray}
  \left\| \phi_{n,j}^{\prime }-\phi_{n,j} \right\|_{H} & \leq & a_j \left\| C_n - C \right\|_{\mathcal{L} (H)} \leq 2 \sqrt{2} \Lambda_{k_n}  \left\| C_n - C \right\|_{\mathcal{S} (H)},
  \end{eqnarray}
  which implies that
  \begin{eqnarray}
\mathcal{P} \left( \sup_{1\leq j\leq k_{n}}\| \phi_{n,j}^{\prime }-\phi_{n,j} \|_{H} \geq \eta \right) & \leq& \mathcal{P} \left( \left\| C_n - C \right\|_{\mathcal{S} (H)} \geq  \frac{\eta}{2 \sqrt{2}  \Lambda_{k_n}} \right). \label{32AA}
\end{eqnarray}

Thus, since $\left\| X_0 \right\|_{H}$ is bounded, from Theorem 4.2 on pp. 99--100 of \cite{Bosq00}, under \textbf{Assumption A1} for any $\eta > 0$, and $\beta  > 1/2,$
\begin{eqnarray}
& & \mathcal{P} \left( \frac{n^{1/4}}{\left(\ln(n) \right)^{\beta}} \sup_{1\leq j\leq k_{n}}\| \phi_{n,j}^{\prime }-\phi_{n,j} \|_{H} \geq \eta \right)
\nonumber \\
&
\leq &\mathcal{P} \left(\left\| C_n - C \right\|_{\mathcal{S} (H)} \geq  \frac{\eta}{2 \sqrt{2} \Lambda_{k_n}}  \frac{\left(\ln(n) \right)^{\beta}}{n^{1/4}}  \right) \nonumber \\
& \leq &  4 \exp \left(- \frac{n \dfrac{\eta^2}{8 \Lambda_{k_n}^{2}} \dfrac{\left( \ln(n) \right)^{2\beta}}{n^{1/2}}}{\gamma_1 + \delta_1 \dfrac{\eta}{2 \sqrt{2} \Lambda_{k_n}}  \dfrac{\left(\ln(n) \right)^{\beta}}{n^{1/4}} } \right)  \nonumber \\
&= &  \mathcal{O}\left(n^{-\frac{\eta ^{2}}{\gamma_1+\eta \delta_1\left(\frac{\ln(n)}{n}\right)^{1/2}}}\right),\quad n\rightarrow \infty.
\label{fi}
\end{eqnarray}
\noindent Thus,  taking  $\eta^{2}> \gamma_1+\delta_1 \eta,$ sequence (\ref{fi}) is summable, and  applying Borel-Cantelli Lemma we arrive to the desired result.\hfill \hfill $\blacksquare$
\end{proof}

\begin{proposition}
\label{prr1}
Under the conditions of Lemma \ref{lemmanew},
considering \textbf{Assumption A1},     for
$\beta > \frac{1}{2},$ and $n$ sufficiently large,
\begin{equation}
\frac{n^{1/4}}{\left(\ln(n) \right)^{\beta}}\sup_{j \geq 1} \left| D_{n,j}- D_j \right| \to^{a.s.} 0,\quad n\rightarrow \infty,
\label{resscc}
\end{equation}
 \noindent where $\left\lbrace D_{n,j}, \ j \geq 1 \right\rbrace$  are defined in equation (\ref{38aaa}), and $\left\lbrace D_{j}, \ j \geq 1 \right\rbrace$ are given by $D_{j}=D_{X}(\phi_{j})(\phi_{j}),$ for every  $j\geq 1,$ with, as before,  $\{\phi_{j},\ j\geq 1\}$ being  the system of eigenvectors of $C_{X}.$
 \end{proposition}
 \begin{remark}
 \label{remratec}
 Note that, under the conditions of Lemma \ref{lemmanew},  since $\| X_{0}\|_{H}$ is bounded, from Theorem \ref{theorem2}, the rate of convergence to zero in Lemma \ref{lem1} and Proposition \ref{prr1} can be improved up to the value $\frac{\left(\ln(n) \right)^{\beta}}{n^{\gamma}},$ for $\beta >1/2,$
  and $\gamma \in(1/4, 1/2),$ allowing  larger values of the truncation order $k_{n},$ given by (\ref{trunc}), for a fixed sample size $n.$ \end{remark}
\begin{proof}
From Theorem \ref{theorem2}, there exists an $n_{0}$ such that for $n\geq n_{0},$ $D_{n}$  is a  Hilbert-Schmidt operator. Then, for  $n\geq n_{0},$ and  for every $j,$ 
\begin{eqnarray}
 &&\frac{n^{1/4}}{\left(\ln(n) \right)^{\beta}}\left| D_{n,j} - D_j \right|=\frac{n^{1/4}}{\left(\ln(n) \right)^{\beta}}\nonumber\\
& &\hspace*{1cm}\times
  \left|D_{n}(\phi_{n,j})(\phi_{n,j})-D_{n}(\phi_{n,j})(\phi_{j})
+D_{n}(\phi_{n,j})(\phi_{j})\right.\nonumber\\ &&\left.\hspace*{2cm}-D(\phi_{n,j})(\phi_{j})+D(\phi_{n,j})(\phi_{j})
-D(\phi_{j})(\phi_{j})\right|\nonumber\\
& &\leq \frac{n^{1/4}}{\left(\ln(n) \right)^{\beta}}\left[\|D_{n}(\phi_{n,j})\|_{H}\|\phi_{n,j}-\phi_{j}\|_{H}+\|(D_{n}-D)(\phi_{n,j})\|_{H}\|\phi_{j}\|_{H}\right.
\nonumber\\
& &
\left.\hspace*{1cm}+\|D(\phi_{n,j}-\phi_{j})\|_{H}\|\phi_{j}\|_{H}\right]\nonumber\\
& &\leq \frac{n^{1/4}}{\left(\ln(n) \right)^{\beta}}\left[\|D_{n}\|_{\mathcal{L}(H)}\|\phi_{n,j}-\phi_{j}\|_{H}+\|D_{n}-D\|_{\mathcal{L}(H)}\right.\nonumber\\
& &
\left.\hspace*{1cm}+\|D\|_{\mathcal{L}(H)}\|\phi_{n,j}-\phi_{j}\|_{H}\right],
\label{dd2bb}
\end{eqnarray}
\noindent where, as before, $\{\phi_{j}, \ j \geq 1\}$ and $\{\phi_{n,j},\ j \geq 1\}$  respectively denote the theoretical  and empirical eigenvectors of $C_{X},$   and $\|\cdot\|_{\mathcal{L}(H)}$ denotes the norm in the space of bounded linear operators.

From Theorem \ref{theorem2},
\begin{equation}
\frac{n^{1/4}}{\left(\ln(n) \right)^{\beta}}\|D_{n}-D\|_{\mathcal{L}(H)}\leq \frac{n^{1/4}}{\left(\ln(n) \right)^{\beta}}
\|D_{n}-D\|_{\mathcal{S}(H)}\to^{a.s.} 0,
\label{sc}
\end{equation}

\noindent and,  for $n$ sufficiently large, $\|D_{n}\|_{\mathcal{L}(H)} <\infty.$ Furthermore,
from Lemma \ref{lemmanew} (see equation (\ref{leqkk})),
\begin{equation}
\frac{n^{1/4}}{\left(\ln(n) \right)^{\beta}}\sup_{1\leq j\leq k_{n}}\|\phi_{n,j}-\phi_{j}\|_{H}\to^{a.s.} 0.
\label{eiv}
\end{equation}
\noindent Hence, from equations (\ref{sc}) and (\ref{eiv}), taking the supremum in $j$ at the left-hand side of equation (\ref{dd2bb}),
we obtain equation (\ref{resscc}).
\hfill \hfill $\blacksquare$

\end{proof}

\section{Strong-consistency  in the  trace operator norm}
\label{sec:4}

In the subsequent developments, we assume that 
the eigenvectors \linebreak $\{\phi_{j},\ j\geq 1\}$ of $C_{X}$ are  unknown. The diagonal componentwise estimator of $\rho $ formulated below is then defined in terms of the empirical eigenvectors.  

The following  condition is assumed  in the remainder of this section:

\vspace{0.4cm}

\noindent \textbf{Assumption A2.} The empirical eigenvalue $C_{n,k_n} > 0~a.s,$ where $k_n$ is the truncation parameter satisfying the conditions established in (\ref{trunc}).

\vspace{0.5cm}

Under \textbf{Assumption A2}, from a functional sample of size $n,$ $X_{0},\dots, X_{n-1},$ the following estimator $\widehat{\rho}_{k_{n}}$ of $\rho$ is formulated:

\begin{equation}
\widehat{\rho}_{k_n} = \displaystyle \sum_{j=1}^{k_n} \rho_{n,j} \phi_{n,j} \otimes \phi_{n,j} = \displaystyle \sum_{j=1}^{k_n} \frac{D_{n,j}}{C_{n,j}} \phi_{n,j} \otimes \phi_{n,j}, \label{140}
\end{equation}
\noindent where, for each $j \geq 1$ and $n \geq 2$,
\begin{equation}
\rho_{n,j} =  \frac{D_{n,j}}{C_{n,j}}  = \frac{\frac{1}{n-1} \displaystyle \sum_{i=0}^{n-2} \left\langle X_{i},\phi_{n,j}\right\rangle_{H}\left\langle X_{i+1},\phi_{n,j}\right\rangle_{H}}{\frac{1}{n} \displaystyle \sum_{i=0}^{n-1} \left[\left\langle X_{i},\phi_{n,j}\right\rangle_{H}\right]^{2}} = \frac{n}{n-1} \frac{\displaystyle \sum_{i=0}^{n-2} X_{i,j,n} X_{i+1,j,n}}{\displaystyle \sum_{i=0}^{n-1} X_{i,j,n}^{2}}. \label{141}
\end{equation}
   The strong-consistency of $\widehat{\rho}_{k_{n}},$ in the trace norm, is derived in the following result under suitable conditions.  
\begin{theorem}
\label{proposition5}
Under \textbf{Assumption A2}, and  the conditions assumed in  Lemma \ref{lemmanew}, if $\rho$ is a positive trace operator,  then,  $\widehat{\rho}_{k_n},$ introduced in equations (\ref{140})--(\ref{141}), is strongly-consistent in the trace norm, i.e.,

\begin{equation}\|\rho-\widehat{\rho}_{k_{n}}\|_{1}\to^{a.s.} 0,\quad n\rightarrow \infty,\label{mreq1}
\end{equation}

\noindent  for 
$k_{n}=o\left(\frac{n^{1/4}}{\left(\ln(n) \right)^{\beta}}\right),$ as $n \rightarrow \infty.$ Here, $\|\cdot\|_{1}$ denotes the trace operator norm. Consequently, $\widehat{\rho}_{k_{n}}$ is also strongly-consistent in the spaces $\mathcal{S}(H)$ and $\mathcal{L}(H).$ 
\end{theorem}

\begin{proof}
From Theorem \ref{theorem2}, we have, for $n$ sufficiently large,
\begin{eqnarray}&&\|D_{n}C_{n}^{-1}-D_{X}C_{X}^{-1}\|_{\mathcal{S}(H)}\nonumber\\
& &=
\| D_{n}C_{n}^{-1}-D_{X}C_{n}^{-1}+D_{X} C_{n}^{-1}-D_{X}C_{X}^{-1}\|_{\mathcal{S}(H)}\nonumber\\
& & \leq\| D_{n}C_{n}^{-1}-D_{X}C_{n}^{-1}\|_{\mathcal{S}(H)}+
\|D_{X}C_{n}^{-1}-D_{X}C_{X}^{-1}\|_{\mathcal{S}(H)}\nonumber\\
&&= \|(D_{n}-D_{X})C_{n}^{-1}\|_{\mathcal{S}(H)}+
\|D_{X}(C_{n}^{-1}-C_{X}^{-1})\|_{\mathcal{S}(H)},\nonumber\\
\label{eqproof1}
\end{eqnarray}
\noindent since from such a theorem,  $D_{X} C_{n}^{-1}\in \mathcal{S}(H)$ for $n$ sufficiently large. Again, from Theorem \ref{theorem2}, equation (\ref{eqproof1}) tends to zero as $n\rightarrow \infty,$ a.s., which means that $\|D_{n}C_{n}^{-1}-D_{X}C_{X}^{-1}\|_{\mathcal{S}(H)}$ also converges to zero a.s. From  Theorem \ref{theorem2} and equation (\ref{eqproof1}), we also have that    there exists an $n_{0}\in \mathbb{N}$ such that,   for $n\geq n_{0},$  $D_{n}C_{n}^{-1}$ is a positive trace operator almost surely. Considering now  $n\geq n_{0},$ for any orthonormal basis $\{\varphi_{j}, \ j \geq 1 \}$  on $H,$

\begin{eqnarray}\|D_{n}C_{n}^{-1}\|_{1}&=&\sum_{j=1}^{\infty}\left\langle D_{n}C_{n}^{-1}(\varphi_{j}),\varphi_{j}\right\rangle_{H}.\label{eqproof3}
\end{eqnarray}
\noindent Note also that \begin{eqnarray}
\|D_{X}C_{X}^{-1}\|_{1}&=&\sum_{j=1}^{\infty}\left\langle D_{X}C_{X}^{-1}(\varphi_{j}),\varphi_{j}\right\rangle_{H}.\nonumber\end{eqnarray}

Furthermore,  from Lemma \ref{lem1}, as $n\rightarrow \infty,$ $C_{n,j}^{-1}$ converges a.s. to $C_{j}^{-1},$ uniformly in $j,$ with convergence rate at least $\frac{\left(\ln(n) \right)^{\beta}}{n^{1/4}},$  for $\beta >1/2.$  The same assertion holds for the  uniform a.s. convergence of $D_{n,j}$ to $D_{j},$  as $n\rightarrow \infty,$  from Proposition \ref{prr1}. Consequently, for every $j,$ as $n\rightarrow \infty,$
$D_{n,j}C_{n,j}^{-1}$ converges a.s.  to $D_{j}C_{j}^{-1},$  with  uniform a.s.  rate of convergence at least $\frac{\left(\ln(n) \right)^{\beta}}{n^{1/4}},$  for $\beta >1/2.$ 
Equivalently,  from equation (\ref{c11b}) in Lemma \ref{lem1},   and equation (\ref{resscc})  in Proposition \ref{prr1}, 
\begin{eqnarray}\|\widehat{\rho}_{k_{n}}\|_{1} &=&\sum_{j=1}^{k_{n}}D_{n,j}C_{n,j}^{-1}\nonumber\\ 
 & & \leq f(n)= \mathcal{O}\left(k_{n}\frac{\left(\ln(n) \right)^{\beta}}{n^{1/4}}+\sum_{j=1}^{\infty}D_{j}C_{j}^{-1}\right),\ n\rightarrow \infty \ \mbox{a.s}.
\label{eqproof5}
\end{eqnarray}
   Since    $k_{n}=o\left(\frac{n^{1/4}}{\left(\ln(n) \right)^{\beta}}\right),$    equation (\ref{eqproof5}) converges  a.s., when $n \to \infty$,
to \begin{equation}\|D_{X}C_{X}^{-1}\|_{1}=\sum_{j=1}^{\infty}D_{j}C_{j}^{-1}=\sum_{j=1}^{\infty}D_{X}(\phi_{j})(\phi_{j})[C_{X}(\phi_{j})(\phi_{j})]^{-1}=\|\rho \|_{1}.\label{eqlimas}
\end{equation}  

For two positive trace operators, $\mathcal{K}$ and $\mathcal{T},$ the following identities are satisfied, for any orthonormal basis  $\{\varphi_{j}, \ j \geq 1 \}$ on $H,$

\begin{eqnarray} \|\mathcal{K}-\mathcal{T}\|_{1}&=&\mbox{trace}(|\mathcal{K}-\mathcal{T}|) = \max(\|\mathcal{K}\|_{1}, \|\mathcal{T}\|_{1})-\min(\|\mathcal{K}\|_{1}, \|\mathcal{T}\|_{1})\nonumber\\
 \|\mathcal{K}\|_{1}&=&\sum_{j=1}^{\infty}\left\langle \mathcal{K}(\varphi_{j}),
\varphi_{j}\right\rangle_{H},\quad  \|\mathcal{T}\|_{1}=\sum_{j=1}^{\infty}\left\langle \mathcal{T}(\varphi_{j}),
\varphi_{j}\right\rangle_{H}.\nonumber
\end{eqnarray}
\noindent In particular, for $n\geq n_{0},$    \begin{eqnarray}& &\|\rho -\widehat{\rho}_{k_{n}}\|_{1}=\max\left(\|\rho\|_{1},\|\widehat{\rho}_{k_{n}}\|_{1}\right) -\min\left(\|\rho\|_{1},\|\widehat{\rho}_{k_{n}}\|_{1}\right).\label{feproof}\end{eqnarray}
The limit (\ref{mreq1}) then follows from  equations (\ref{eqproof5})--(\ref{feproof}).

\hfill \hfill $\blacksquare$
\end{proof}

The strong consistency in $H$ of the associated ARH(1) plug-in predictor $\widehat{\rho}_{k_n}(X_{n-1})$
of $X_{n}$ is now derived.

\begin{corollary}
Under the conditions of Theorem \ref{proposition5}
\label{proposition6}
\begin{equation}
\|\widehat{\rho}_{k_n}(X_{n-1})-\rho(X_{n-1})\|_{H} \to^{a.s.} 0,\quad n\rightarrow \infty.\label{eqfr}
\end{equation}

\end{corollary}

\begin{proof}
The proof directly  follows from Theorem \ref{proposition5}, keeping in mind that
the convergence in the trace norm implies the convergence in the space $\mathcal{L}(H)$ of bounded linear operators. Moreover,  
$$\|\widehat{\rho}_{k_n}(X_{n-1})-\rho(X_{n-1})\|_{H}\leq \|\widehat{\rho}_{k_n}-\rho\|_{\mathcal{L}(H)}\|X_{n-1}\|_{H}\to^{a.s.} 0, \ n\to \infty,$$
\noindent since $\|X_{n-1}\|_{H}<\infty.$

\hfill \hfill $\blacksquare$
\end{proof}

\begin{remark}
Note that, when $\rho $ is Hilbert-Schmidt, but it is not positive trace operator, under the conditions of   Lemma \ref{lemmanew}, and  if \textbf{Assumption A2} holds, the following a.s. inequality is satisfied,  as $n \to
\infty,$ \begin{equation}\|\widehat{\rho}_{k_{n}}-\rho\|_{\mathcal{S}(H)}^{2}\leq 
\|\rho\|_{\mathcal{S}(H)}^{2}-\sum_{j=1}^{\infty}[\rho(\phi_{j})(\phi_{j})]^{2}=\sum_{j\neq k}^{\infty}\left[\frac{D_{X}(\phi_{j})(\phi_{k})}{C_{j}}\right]^{2}<\infty.\label{mrdef}\end{equation}\end{remark}

\section*{Acknowledgments}

This work has been supported in part by project MTM2015--71839--P (co-funded by Feder funds),
of the DGI, MINECO, Spain.

\end{document}